\documentclass[11pt]{article}
\usepackage{amsthm,amsfonts,amsmath,amssymb}
\usepackage{mdwlist}
\newtheorem{thm}{Theorem}[section]

\newtheorem{lem}[thm]{Lemma}

\newtheorem{cons}[thm]{Construction}

\def\demo{\noindent{\bf Proof}\hskip10pt}

\def\intl{[\hskip-5pt[\hskip 2pt}
\def\intr{\hskip 2pt ]\hskip-5pt]\hskip 2pt}

\def\qed{\hfill $\Box$}
\def\lg{\langle}
\def\rg{\rangle}
\def\rr#1{\item[{\rm (#1)}]}

\def\CD{\mathcal{CD}}
\def\L{\mathcal{L}}

\def\span{\hbox{\rm span}}
\begin{document}
\title{Groups whose Chermak-Delgado lattice is a subgroup lattice of an elementary abelian $p$-group\thanks{This work was supported by NSFC (No. 11971280 \&11771258)}}
\author{Lijian An\\
Department of Mathematics,
Shanxi Normal University\\
Linfen, Shanxi 041004, P. R. China\\
 }

\maketitle

\begin{abstract}

The Chermak-Delgado lattice of a finite group $G$ is a self-dual sublattice of the subgroup lattice of $G$.
In this paper, we focus on finite groups whose Chermak-Delgado lattice is a subgroup lattice of an elementary abelian $p$-group. We prove that such groups are nilpotent of class $2$. We also prove that, for any elementary abelian $p$-group $E$, there exists a finite group $G$ such that the Chermak-Delgado lattice of $G$ is a subgroup lattice of $E$.

\medskip


\end{abstract}

\baselineskip=16pt

\section{Introduction}

 Suppose that $G$ is a finite group, and $H$ is a subgroup of $G$. The Chermak-Delgado measure of $H$ (in $G$) is denoted by $m_G(H)$, and defined as $m_G(H)=|H|\cdot |C_G(H)|.$
The maximal Chermak-Delgado measure of $G$ is denoted by $m^*(G)$, and defined as $$m^*(G)=\max\{ m_G(H)\mid H\le G\}.$$  Let $$\mathcal{CD}(G)=\{ H\mid m_G(H)=m^*(G)\}.$$ Then the set $\mathcal{CD}(G)$ forms a sublattice of the subgroup lattice of $G$, which is
called the Chermak-Delgado lattice of $G$.
It was first introduced by Chermak and Delgado \cite{CD}, and revisited by Isaacs \cite{I}.  In the last years, there has been a growing interest in understanding this lattice (see e.g. [1-6], [8-9], [12-18]).
%
%
%
%

Notice that a Chermak-Delgado lattice is always self-dual. It is natural to ask the question: which types of self-dual lattices can be as Chermak-Delgado lattices of finite groups. Some special cases of this question are proposed and solved. In \cite{BHW1}, it is proved that, for any integer $n$, a chain of length $n$ can be a Chermak-Delgado lattice of a finite $p$-group.

A quasi-antichain is a lattice consisting of a maximum, a minimum, and the atoms of the lattice. The width of a quasi-antichain is the number of atoms.
For a positive integer $w\ge 3$, a quasi-antichain of width $w$ is denoted by $\mathcal{M}_{w}$. In \cite{BHW2}, it was proved that $\mathcal{M}_{w}$ can be as a Chermak-Delgado lattice of a finite group if and only if $w=1+p^a$ for some positive integer $a$ and some prime $p$.
The following theorem gives more self-dual lattices which can be as Chermak-Delgado lattices of finite groups.

 \begin{thm}{\rm (\cite{ABQW})}\label{th=abqw}
 If $\mathcal{L}$ is a Chermak-Delgado lattice of a finite $p$-group $G$ such that both $G/Z(G)$ and $G'$ are elementary abelian, then are $\mathcal{L}^+$ and $\mathcal{L}^{++}$, where $\mathcal{L}^+$ is a mixed $3$-string with center component isomorphic to $\mathcal{L}$ and the remaining components being $m$-diamonds {\rm (}a lattice with subgroups in the configuration of an $m$-dimensional cube{\rm )}, $\mathcal{L}^{++}$ is a mixed $3$-string with center component isomorphic to $\mathcal{L}$ and the remaining components being lattice isomorphic to $\mathcal{M}_{p+1}$.
\end{thm}

%
%


For a finite group $G$, we use $\L(G)$ to denoted the subgroup lattice of $G$. We use $E_{p^n}$ to denote the elementary abelian $p$-group of order $p^n$. It is well-known that $\L(E_{p^n})$ is self-dual. Let $G$ be an extra-special $p$-group of order $p^{2n+1}$. Then $\CD(G)$ is isomorphic to $\L(E_{p^{2n}})$ (see \cite[Example 2.8]{GL}).
In this paper, we focus on finite groups whose Chermak-Delgado lattice is isomorphic to $\L(E_{p^n})$. The main results are:

\medskip

\noindent {\bf Theorem A.}  Let $G$ be a finite group with $G\in\mathcal{CD}(G)$. Suppose that $\mathcal{CD}(G)$ is isomorphic to $\L(E_{p^n})$, where $n\ge 2$. Then $G=P\times Q$, where $P$ is the Sylow $p$-subgroup of $G$ such that $P/Z(P)$ is elementary abelian, $Q$ is the abelian Hall $p'$-subgroup of $G$. Moreover, $\mathcal{CD}(G)\cong \mathcal{CD}(P)$ as lattice.

\medskip

\noindent {\bf Theorem B.} For any integer $n$ and a prime $p$, there exists a special $p$-group $G$ such that $\CD(G)$ is isomorphic to $\L(E_{p^n})$.

\medskip

For a Chermak-Delgado lattice, the following properties is basic and is often used in this paper. We will not point out when we use them.

\begin{thm}{\rm \cite{CD}}\label{basic} Suppose that $G$ is a finite group and $H,K\in\CD(G)$.
\begin{itemize*}
  \rr{1} $\lg H,K\rg=HK$. Hence a Chermak-Delgado lattice is modular.
  \rr{2} $C_G(H\cap K)=C_G(H)C_G(K)$.
   \rr{3} $C_G(H)\in\mathcal{CD}(G)$ and $C_G(C_G(H))=H$. Hence a Chermak-Delgado lattice is self-dual.
  \rr{4} Let $M$ be the maximal member of $\mathcal{CD}(G)$. Then $M$ is characteristic in $G$ and $\mathcal{CD}(M)=\mathcal{CD}(G)$.
  \rr{5} The minimal member of $\mathcal{CD}(G)$ is characteristic, abelian, and contains $Z(G)$.
\end{itemize*}
\end{thm}

\section{Quasi-antichain intervals in Chermak-Delgado lattices}

If $n\ge 2$, then every interval of length $2$ in $\L(E_{p^n})$ is a quasi-antichain of width $p+1$. Hence we start our argument from investigating quasi-antichain intervals in Chermak-Delgado lattices.
Following \cite{BHW2}, we use $\intl L,H\intr$ to denote the interval from $L$ to $H$ in $\mathcal{CD}(G)$.

%

\begin{lem}{\rm (\cite[Proposition 2 \& Theorem 4]{BHW2})}\label{prop3}
Let $G$ be a finite group with an interval $\intl L,H\intr\cong \mathcal{M}_w$ in $\mathcal{CD}(G)$, where $w\ge 3$, and $K$ be an atom of the quasi-antichain. Then $K\trianglelefteq H$, $L\trianglelefteq H$, and there exists a prime $p$ and positive
integers $a, b$ with $b\le a$ such that $H/L$ is elementary
abelian $p$-groups of order $p^{2a}$, $|H/K|=|K/L|=p^a$ and $w = p^b + 1$.
\end{lem}

\begin{lem}\label{gengeralize}
Let $G$ be a finite group with an interval $\intl L,H\intr$ of length $l$ in $\mathcal{CD}(G)$. Suppose that every interval of length $2$ in $\intl L,H\intr$ is a quasi-antichain of width $\ge 3$. Then $L\trianglelefteq H$, and there exists a prime $p$ and a positive
integer $a$ such that $H/L$ is elementary
abelian $p$-groups of order $p^{al}$. Moreover, if $\intl L_1,H_1\intr$ is an interval of length $2$ in $\intl L,H\intr$, then the width of $\intl L_1,H_1\intr$ is $p^b+1$ for some integer $b$.
\end{lem}
\demo We use $l_{H/L}$ to denote the length of $\intl L,H\intr$. If $l=l_{H/L}=2$, then the conclusions follow from Lemma \ref{prop3}. In the following, we may assume that $l \ge 3$. Let $L=J_0<J_1<J_2<\dots<J_l=H$ be a maximal chain.

For $0\le i\le l-2$, $\intl J_i,J_{i+2}\intr$ is a quasi-antichain of width $\ge 3$. By Lemma \ref{prop3}, there exist  primes $p_i$ and positive
integers $a_i$ such that, $|J_{i+2}/J_{i+1}|=|J_{i+1}/J_i|=p_i^{a_i}$. It is easy to see that these $p_i$ are coincide to a prime $p$ and these $a_i$ are coincide to an integer $a$. Hence $|H:L|=p^{al}$.
   
Since $\intl L,J_2\intr$ is a quasi-antichain of width $\ge 3$, there is an atom $L<K_1<J_2$ such that $K_1\ne J_1$.
Since $l_{H/J_1}=l_{H/K_1}=l-1$, by induction, $J_1\trianglelefteq H$, $K_1\trianglelefteq H$ and $H/J_1$ and $H/K_1$ are elementary $p$-groups. It follows that $L=J_1\cap K_1 \unlhd H$ and $H/L\lesssim H/J_1\times H/K_1$ is an elementary abelian $p$-group of order $p^{al}$.

 If $\intl L_1,H_1\intr$ is an interval of length $2$ in $\intl L,H\intr$, then, by Lemma \ref{prop3}, the width of $\intl L_1,H_1\intr$ is $p^b+1$ for some integer $b$.
\qed

\begin{thm}\label{nnn}
Let $G$ be a finite group with $G\in\mathcal{CD}(G)$. Suppose that every interval of length $2$ in $\mathcal{CD}(G)$ is a quasi-antichain of width $\ge 3$. Then there exists a prime $p$ such that the width of an interval of length $2$ in $\CD(G)$ is $p^b+1$ for some integer $b$.
 Moreover, $G=P\times Q$, where $P$ is the Sylow $p$-subgroup of $G$, nilpotent of class $2$, $Q$ is the abelian Hall $p'$-subgroup of $G$, and $\mathcal{CD}(G)\cong \mathcal{CD}(P)$ as lattice.
\end{thm}
\demo By Lemma \ref{gengeralize}, $G/Z(G)$ ia an elementary abelian $p$-group for some prime $p$ and the width of an interval of length $2$ in $\CD(G)$ is $p^b+1$ for some integer $b$. Hence $G$ is nilpotent. Let $P$ be the Sylow $p$-subgroup and $Q$ the Hall $p'$-subgroup. Then $G=P\times Q$. Obviously, $Q\in Z(G)$. Therefore $\mathcal{CD}(G)\cong \mathcal{CD}(P)\times \mathcal{CD}(Q)\cong \mathcal{CD}(P)$.
\qed

\medskip

\noindent{\bf Proof of Theorem A.} Since $\mathcal{CD}(G)$ is isomorphic to $\L(E_{p^n})$, where $n\ge 2$, every interval of length $2$ in $\mathcal{CD}(G)$ is a quasi-antichain of width $p+1$. By Theorem \ref{nnn}, there exists a prime $p_1$ such that the width of an interval of length $2$ in $\CD(G)$ is $p_1^b+1$ for some integer $b$. Hence $p_1=p$ and $b=1$. Other results hold obviously.
\qed

\section{Proof of Theorem B}

\begin{lem}\rm{\cite[Lemma 1.43]{I}}\label{ii}
Suppose that $G$ is finite group. If $H,K\le G$, then $$m_G(H)\cdot m_G(K) \le m_G(\lg H,K\rg)\cdot m_G(H\cap K).$$ Moreover, equality occurs if and only if $\lg H,K\rg=HK$ and $C_G(H\cap K)=C_G(H)C_G(K)$.
\end{lem}

\begin{lem}\label{lem relation}
Suppose that $G$ is finite group. If $K\le H\le G $, then $$\frac{m_H(K)}{m_G(K)}\le \frac{m_H(H)}{m_G(H)}. $$ Moreover, equality occurs if and only if $C_G(K)\le HC_G(H)$.
\end{lem}
\demo By calculation,
$$ \frac{m_H(K)}{m_G(K)}=\frac{|K|\cdot |C_H(K)|}{|K|\cdot |C_G(K)|}=\frac{|C_H(K)|}{|C_G(K)|}=\frac{|H\cap C_G(K)|}{|C_G(K)|} =\frac{|H|}{|HC_G(K)|}$$
and
$$\frac{m_H(H)}{m_G(H)}=\frac{|H|\cdot |C_H(H)|}{|H|\cdot |C_G(H)|}=\frac{|C_H(H)|}{|C_G(H)|}=\frac{|H\cap C_G(H)|}{|C_G(H)|} =\frac{|H|}{|HC_G(H)|}.$$
\smallskip

\noindent
Since $K\le H$, $C_G(H)\le C_G(K)$. Hence $|HC_G(H)|\le |HC_G(K)|$, where equality occurs if and only if $C_G(K)\le HC_G(H)$. Thus
$$ \frac{m_H(K)}{m_G(K)}=\frac{m_H(H)}{m_G(H)}\frac{|HC_G(H)|}{|HC_G(K)|}\le \frac{m_H(H)}{m_G(H)},$$ where equality occurs if and only if $C_G(K)\le HC_G(H)$.
\qed

\begin{lem}\label{maximal element}
Suppose that $G$ is finite group, $H\le G $ such that $G=HC_G(H)$. If $H\in\CD(H)$, then $H$ contains in the unique maximal member of $\CD(G)$.
\end{lem}
\demo
Let $M$ be the unique maximal member of $\CD(G)$. Since $C_G(M\cap H)\le G=HC_G(H)$, by Lemma \ref{lem relation}, $$\frac{m_H(M\cap H)}{m_G(M\cap H)}=\frac{m_H(H)}{m_G(H)}.$$
Since $H\in\CD(H)$, $m_H(M\cap H)\le m_H(H)$. It follows that
$m_G(M\cap H)\le m_G(H).$ By Lemma \ref{ii},
$$m_G(H)\cdot m_G(M) \le m_G(\lg H,M\rg)\cdot m_G(M\cap H).$$ It follows that $m_G(\lg H,M\rg)\ge m_G(M)=m^*(G)$. Hence $\lg H,M\rg\in \CD(G)$. Since $M$ is maximal, $H\le M$.
\qed

\begin{thm}\label{thm sub}
Suppose $H\in \CD(G)$. Then $\CD(H)$ is just the interval $\intl Z(H), H\intr$ in $\CD(G)$.
\end{thm}
\demo Since $H\in \CD(G)$, $m_G(H)=m^*(G)$. By Lemma \ref{lem relation},
\begin{equation}\label{eq}
  m_H(H)\ge m_G(H)\frac{m_H(K)}{m_G(K)}=m^*(G)\frac{m_H(K)}{m_G(K)}\ge m_H(K)
\end{equation}
 for all $K\le H$.
It follows that $H\in\CD(H)$ and $m^*(H)=m_H(H)$. Moreover, ``=" holds in Equation (\ref{eq}) if and only if $C_G(K)\le HC_G(H)$ and $m_G(K)=m^*(G)$.
Notice that $C_G(K)\le HC_G(H)$ if and only if $K\ge C_G(HC_G(H))=H\cap C_G(H)=Z(H)$, and $m_G(K)=m^*(G)$ if and only if $K\in\CD(G)$. $K\in\CD(H)$ if and only if
\begin{center}
``=" holds in Equation (\ref{eq})
\end{center}
if and only if $K\in \intl Z(H), H\intr$. Hence $\CD(H)$ is just the interval $\intl Z(H), H\intr$.
\qed

\begin{lem}\label{lemma 1} Let $\mathbb{F}$ be a field and $A,B$ be $n\times n$ matrices, where $n\ge 3$. If $AZ=ZB$ for all anti-symmetric matrix $Z$, then $A, B$ are scalar matrices.
\end{lem}
\demo Let $A=(a_{ij})$ and $B=(b_{ij})$. Since $A(E_{ij}-E_{ji})=(E_{ij}-E_{ji})B$ for $ i\ne j$, by calculation, we have $a_{ii}=b_{jj}$ and
\begin{center}
$a_{ki}=a_{kj}=b_{ik}=b_{jk}=0$, if $k\ne i,j.$
\end{center}
Since $n\ge 3$, all $a_{ki}$ and $b_{ki}$ are zero if $k\ne i$, and all $a_{ii}$ and $b_{ii}$ are coincide. Hence  $A$ and $ B$ are scalar matrices.
\qed

\begin{cons}\label{cons} For a prime $p$, let $P=\lg x,y,w; z_1,z_2,z_3\mid x^p=y^p=w^p=1, [x,y]=z_1,[y,w]=z_2,[w,x]=z_3, z_i^p=[z_i,x]=[z_i,y]=[z_i,w]=1\ \mbox{where}\ i=1,2,3 \rg$. Then it is easy to check that $\Phi(P)=Z(P)=P'=\lg z_1,z_2,z_3\rg$ is of order $p^3$, $|P|=p^6$, and $\CD(P)=\{ P,Z(P)\}$.

Let $G_n$ be the group which is the central product of $n$ copies of $P$. Thus $G_n$ has order $p^{3n+3}$ and is generated by $3n$ elements of order $p$, $x_1,\dots,x_n,y_1,\dots,y_n,w_1,\dots,w_n$ subject to the defining relations:
\begin{center}
$[x_i,x_j]=[y_i,y_j]=[w_i,w_j]=[x_i,y_j]=[y_i,w_j]=[w_i,x_j]=1$ if $i\ne j$,

$[x_i,y_i]=z_1$, $[y_i,z_i]=z_2$, $[w_i,x_i]=z_3$,

$[z_j,x_i]=[z_j,y_i]=[z_j,w_i]=1$.
\end{center}
It is easy to Check that $G_n$ is a special $p$-group with $\Phi(G_n)=Z(G_n)=G_n'=\lg z_1,z_2,z_3\rg$.
\end{cons}

\begin{thm}\label{main}
Let $G=G_n$ which is defined in  Construction {\rm \ref{cons}}. Let $P_i=\lg x_i,y_i,w_i\rg$, where $1\le i\le n$. Then
\begin{itemize*}
 \rr1 $G\in \CD(G)$, $m^*(G)=p^{3n+6}$, and $P_i\in \CD(G)$.
 \rr2 If $H\in \CD(G)$, then $|H|=p^{3m+3}$ for some $0\le m\le n$.
 \rr3 We use $F_p$ to denote the finite field $\mathbb{Z}/p\mathbb{Z}$. For a vector $$v=(s_1,s_2,\dots,s_n)$$ of $F_p^n$, we use $v^\varphi$ to denote the subgroup $\lg \alpha_v,\beta_v,\gamma_v,Z(G)\rg$, where
 $$\alpha_v=\prod_{i=1}^{n}x_i^{s_i},\ \beta_v=\prod_{i=1}^{n}y_i^{s_i},\ \gamma_v=\prod_{i=1}^{n}w_i^{s_i}.$$
 Let $\tilde{v}=(t_1,t_2,\dots,t_n)$. Define an inner product on $F_p^n$ with $\lg v,\tilde{v} \rg=\sum_{i=1}^ns_it_i$. Then $[v^\varphi,{\tilde{v}}^\varphi]=1$ if and only if $\lg v,\tilde{v}\rg=0.$
 \rr4
 Suppose that $U$ is an $m$-dimensional subspace of $F_p^n$. We use $U^\varphi$ to denote the subgroup $\prod_{u\in U}{u}^\varphi$ of $G$. Then $|U^\varphi|=p^{3m+3}$.
 \rr5 $U^\varphi\in\CD(G)$. Moreover, let $$U^\perp=\{ v\in F_p^n\mid \lg u,v\rg=0 \ \mbox{for all}\ u\in U\ \}.$$ Then $(U^\perp)^\varphi=C_G(U^\varphi)$.
 \rr6 If $H\in\CD(G)$, then there exists a subspace $U$ of $F_p^n$ such that $H=U^\varphi$.
\end{itemize*}
\end{thm}
\demo (1) It is easy to see that $C_G(P_i)=\prod_{j\ne i}P_j\cong G_{n-1}$ for $1\le i\le n$.  Hence $m_G(G)=m_G(P_i)=p^{3n+6}$.
Since $\CD(P_i)=\{P_i, Z(G)\}$ and $G=P_iC_G(P_i)$, by Lemma \ref{maximal element}, $P_i$ contains in the unique maximal member of $\CD(G)$. Hence $G$ is the unique maximal member of $\CD(G)$, $m^*(G)=m_G(G)=p^{3n+6}$, and $P_i\in\CD(G)$.

(2) It is trivial for $n=1$. Assume that $n\ge 2$. Let $Q_n=C_G(P_n)=\prod_{i=1}^{n-1}P_i\cong G_{n-1}$. If $P_n\le H$, then $C_G(H)\le Q_n$. By Theorem \ref{thm sub}, $C_G(H)\in\CD(Q_n)$. By induction, $|C_G(H)|=p^{3m+3}$ for some $0\le m\le n-1$. Hence $|H|=m^*(G)/|C_G(H)|=p^{3(n-m)+3}$, where $1\le n-m\le n$. If $P_n\not\le H$, then, by above argument, $|HP_n|=p^{3m+3}$ for some $1\le m\le n$.
By Theorem \ref{thm sub}, $H\cap P_n\in\CD(P_n)$. Since $\CD(P_n)=\{ P_n,Z(G)\}$ and $P_n\not\le H$, $H\cap P_n=Z(G)$.  Hence $$|H|=\frac{|HP_n|\cdot |H\cap P_n|}{|P_n|}=p^{3(m-1)+3}.$$

(3)  By calculation, $$[\alpha_v,\beta_{\tilde{v}}]=[\alpha_{\tilde{v}},\beta_{v}]=z_1^{\lg v,\tilde{v} \rg},$$
$$[\beta_v,\gamma_{\tilde{v}}]=[\beta_{\tilde{v}},\gamma_{{v}}]=z_2^{\lg v,\tilde{v} \rg},$$ $$[\gamma_v,\alpha_{\tilde{v}}]=[\gamma_{\tilde{v}},\alpha_{{v}}]=z_3^{\lg v,\tilde{v} \rg}.$$
Since $[\alpha_v,\alpha_{\tilde{v}}]=[\beta_v,\beta_{\tilde{v}}]=[\gamma_v,\gamma_{\tilde{v}}]=1$, $[v^\varphi,{\tilde{v}}^\varphi]=1$ if and only if $\lg v,\tilde{v} \rg=0.$

\medskip

(4) Let
\begin{center}$A=\lg \alpha_{u}\mid u\in U\rg Z(G),$

$B=\lg \beta_{u}\mid u\in U\rg Z(G)$,

 $C=\lg \gamma_{u}\mid u\in U\rg Z(G).$
\end{center}
Then $|A/Z(G)|=|B/Z(G)|=|C/Z(G)|=p^m$. Since $$U^\varphi/Z(G)=A/Z(G)\times B/Z(G)\times C/Z(G),$$ $|U^\varphi/Z(G)|=p^{3m}$. Hence $|U^\varphi|=p^{3m+3}$.

\medskip

(5) Notice that $U^\perp$ is an $(n-m)$-dimensional subspace of $F_p^n$. By (4), $|({U^\perp})^\varphi|=p^{3(n-m)+3}$. By (3), $[({U^\perp})^\varphi,U^\varphi]=1$. Hence
$$m_G(U^\varphi)=|U^\varphi|\cdot |C_G(U^\varphi)|\geqslant |U^\varphi|\cdot |({U^\perp})^\varphi|=p^{3n+6}=m^*(G).$$ It follows that $U^\varphi\in\CD(G)$ and $C_G(U^\varphi)=({U^\perp})^\varphi.$

\medskip

(6) If $n=1$, then the conclusion is trivial. In the following, we assume that $n\ge 2$.
Let $Q_n=C_G(P_n)=\prod_{i=1}^{n-1}P_i\cong G_{n-1}$.

\medskip

Case 1. $H\le Q_n$.

By Theorem \ref{thm sub}, $H\in\CD(Q_n)$. By induction, there exists a subspace $U$ of $F_p^{n-1}\times \{ 0\}\subset F_p^n$ such that $H=U^\varphi$.

\medskip

Case 2. $P_n\le H$.
 
Notice that $C_G(H)\le C_G(P_n)=Q_n$. By Case 1, there exists a subspace $U$ of $F_p^{n-1}\times \{ 0\}\subset F_p^n$ such that $C_G(H)=U^\varphi$. By (5), $H=C_G(C_G(H))=({U^\perp})^\varphi$.

\medskip

Case 3.  $H\not\le Q_n$ and $P_n\not\le H$.

 By (2), $|H|=p^{3m+3}$ for some $1\le m\le n-1$.
Also by (2), $|HQ_n|=p^{3m'+3}$ for some $n-1\le m'\le n$. Since $HQ_n>Q_n$, $m'=n$ and hence $HQ_n=G$.

Let $H_1=HP_n$, $H_2=HP_n\cap Q_n$ and $H_3=H\cap Q_n$.  Since $H\cap P_n\in \intl Z(P_n),P_n\intr$, by Theorem \ref{thm sub}, $H\cap P_n\in\CD(P_n)=\{ P_n,Z(G)\}$. It follows that $H\cap P_n=Z(G)$. Hence $$|H_1|=|HP_n|=\frac{|H|\cdot |P_n|}{|H\cap P_n|}=p^{3(m+1)+3}.$$
Since $H_1Q_n=HQ_n=G$,
$$|H_2|=|H_1\cap Q_n|=\frac{|H_1|\cdot |Q_n|}{|H_1Q_n|}=p^{3m+3}$$
and
$$|H_3|=|H\cap Q_n|=\frac{|H|\cdot |Q_n|}{|HQ_n|}=p^{3(m-1)+3}.$$
By Case 1, there exist an $m$-dimensional subspace $U_2$ and an $(m-1)$-dimensional subspace $U_3$ of $F_p^{n}\times \{ 0\}\subset F_p^n$ such that $H_2=U_2^\varphi$ and $H_3=U_3^\varphi$.
Since $H_3\le H_2$, $U_3\le U_2$. Hence there exists a vector $u\in F_p^{n}\times \{ 0\}\subset F_p^n$ such that $H_2=H_3 u^\varphi$. Let $H^*=H\cap u^\varphi P_n$.
Then $$H=H\cap HP_n=H\cap H_2P_n=H\cap H_3 u^\varphi P_n=H_3(H\cap u^\varphi P_n)=H_3H^*.$$
Since $H^*P_n=(H\cap u^\varphi P_n)P_n=HP_n\cap u^\varphi P_n=u^\varphi P_n$, we may assume that $H^*=\lg \alpha,\beta,\gamma,Z(G)\rg$ where
$$\alpha=\alpha_u x_n^{a_{11}}y_n^{a_{12}}w_n^{a_{13}},\ \beta=\beta_u x_n^{a_{21}}y_n^{a_{22}}w_n^{a_{23}},\ \gamma=\gamma_u x_n^{a_{31}}y_n^{a_{32}}w_n^{a_{33}}.$$
Since $H_2=H_3u^\varphi$, $$Hu^\varphi=HH_3u^\varphi=HH_2=H(HP_n\cap Q_n)=HP_n\cap HQ_n=HP_n.$$
Hence $$H^*u^\varphi=(H\cap u^\varphi P_n)u^\varphi=Hu^\varphi\cap u^\varphi P_n=u^\varphi P_n.$$
It follows that the matrix $(a_{ij})_{3\times 3}$ is invertible.

Let $K=C_G(H)$.  Then $P_n\not\le K$ and $K\not\le Q_n$. Similarly, there exist a vector $v\in F_p^{n}\times \{ 0\}\subset F_p^n$ such that $K=(K\cap Q_n)(K\cap v^\varphi P_n)$, and we may assume that $K\cap v^\varphi P_n=\lg \tilde{\alpha},\tilde{\beta},\tilde{\gamma},Z(G)\rg$ where
$$\tilde{\alpha}=\alpha_v x_n^{b_{11}}y_n^{b_{12}}w_n^{b_{13}},\ \tilde{\beta}=\beta_u x_n^{b_{21}}y_n^{b_{22}}w_n^{b_{23}},\ \tilde{\gamma}=\gamma_u x_n^{b_{31}}y_n^{b_{32}}w_n^{b_{33}}.$$
The matrix $(b_{ij})_{3\times 3}$ is also invertible.

For convenience, in the following, the operation of the group $G$ is written as addition. Hence
$$\left(
  \begin{array}{c}
    \alpha \\
    \beta \\
    \gamma \\
  \end{array}
\right)=\left(
  \begin{array}{c}
    \alpha_u \\
    \beta_u \\
    \gamma_u \\
  \end{array}
\right)+A \left(
  \begin{array}{c}
    x_n \\
    y_n \\
    z_n \\
  \end{array}
\right)$$
where $A=(a_{ij})_{3\times 3}$, and $(\tilde{\alpha},\tilde{\beta},\tilde{\gamma})=(\alpha_v,\beta_v,\gamma_v)+(x_n,y_n,z_n)B^T$ where $B=(b_{ij})_{3\times 3}.$
Notice that $[u^\varphi,P_n]=[v^\varphi,P_n]=0$.  By calculation, we have
\begin{equation}\label{123}
  \left(
  \begin{array}{ccc}
    [\alpha,\tilde{\alpha}] & [\alpha,\tilde{\beta}] & [\alpha,\tilde{\gamma}] \\
  \mbox{[}\beta,\tilde{\alpha} \mbox{]}& [\beta,\tilde{\beta}] & [\beta,\tilde{\gamma}] \\
  \mbox{[}\gamma,\tilde{\alpha} \mbox{]} & [\gamma,\tilde{\beta}] & [\gamma,\tilde{\gamma}] \\
  \end{array}
\right)=\left(
  \begin{array}{ccc}
    [\alpha_u,\alpha_v] & [\alpha_u,\beta_v] & [\alpha_u,\gamma_v] \\
  \mbox{[}\beta_u,\alpha_v \mbox{]}& [\beta_u,\beta_v] & [\beta_u,\gamma_v] \\
  \mbox{[}\gamma_u,\alpha_v \mbox{]} & [\gamma_u,\beta_v] & [\gamma_u,\gamma_v] \\
  \end{array}
\right)+ AZB^T
\end{equation}
where $$Z=\left(
  \begin{array}{ccc}
    [x_n,x_n] & [x_n,y_n] & [x_n,w_n] \\
  \mbox{[}y_n,x_n \mbox{]}& [y_n,y_n] & [y_n,w_n] \\
  \mbox{[}w_n,x_n \mbox{]} & [w_n,y_n] & [w_n,w_n] \\
  \end{array}
\right)=\left(
  \begin{array}{ccc}
    0 & z_1 & z_2 \\
    -z_1 & 0 & -z_3 \\
    z_3 & -z_2 & 0 \\
  \end{array}
\right).$$
Since $[H^*, K\cap v^\varphi P_n]=0$, the left of Equation (\ref{123}) equals to $O_{3\times 3}$.
Since
$$\left(
  \begin{array}{ccc}
    [\alpha_u,\alpha_v] & [\alpha_u,\beta_v] & [\alpha_u,\gamma_v] \\
   \mbox{[}\beta_u,\alpha_v \mbox{]} & [\beta_u,\beta_v] & [\beta_u,\gamma_v] \\
    \mbox{[}w_n,x_n \mbox{]} & [\gamma_u,\beta_v] & [\gamma_u,\gamma_v] \\
  \end{array}
\right)=\lg u,v\rg Z,$$
we have
\begin{equation}\label{456}
  AZB^T=-\lg u,v\rg Z
\end{equation}
Since $A$ is invertible, we have $ZB^T=-\lg u,v\rg A^{-1}Z$. By Lemma \ref{lemma 1}, $A$ is a scalar matrix. Let $u^*=u+(0,\cdots,0,a_{11})$. Then $H^*=(u^*)^\varphi$. Let $U=U_3+\span(u^*)$. Then $H=H_3 H^*=U^\varphi$.
\qed

\medskip

\noindent{\bf Proof of Theorem B.} Let $G=G_n$ which is defined in  Construction {\rm \ref{cons}}. By Theorem \ref{main}, $\CD(G)$ is isomorphic to the subspace lattice of $F_p^n$, which is isomorphic to $\L(E_{p^n})$. \qed

%

\end{document}